\documentclass[letterpaper,10pt,conference]{ieeeconf}

\IEEEoverridecommandlockouts

\usepackage{amssymb,euscript,latexsym,amsmath,graphicx,bm}
\usepackage[percent]{overpic}
\usepackage{rotating}

\usepackage{algorithm, algorithmic}
 
\newcommand{\bfi}{\bfseries\itshape}

\renewcommand{\epsilon}{\varepsilon}

\makeatletter
\def\thefigure{
\arabic{figure}}
\def\fps@figure{h, t}
\def\theequation{
\arabic{equation}}
\makeatother
\DeclareMathOperator*{\ext}{ext}

\newtheorem{theorem}{Theorem}

\newtheorem{dfn}[theorem]{Definition}

\def\pder#1#2{{\frac {\partial #1} {\partial #2}}}

\begin{document}
\title{\LARGE \bf Controlled Lagrangians and Potential Shaping for \\
Stabilization of Discrete Mechanical Systems}

\author {Anthony M. Bloch
\\ Department of Mathematics
\\ University of Michigan
\\ Ann Arbor, MI 48109
\\ {\tt\small abloch@umich.edu}
\\
\\
Jerrold E. Marsden
\\ Control and Dynamical Systems
\\ California Institute of Technology 107-81
\\ Pasadena, CA 91125
\\ {\tt\small marsden@cds.caltech.edu}
\and 
Melvin Leok
\\ Department of Mathematics
\\ Purdue University
\\ West Lafayette, IN 47907
\\ {\tt\small mleok@math.purdue.edu}
\\
\\
Dmitry V. Zenkov
\\ Department of Mathematics
\\ North Carolina State University \\ Raleigh, NC 27695 \\
{\tt\small dvzenkov@unity.ncsu.edu}
}

\maketitle

\begin{abstract}
\noindent
The method of controlled Lagrangians for discrete mechanical systems is extended to include potential shaping in order to achieve complete state-space asymptotic stabilization. New terms in the controlled shape equation that are necessary for matching in the discrete context are introduced. The theory is illustrated with the problem of stabilization of the cart-pendulum system on an incline. 
We also discuss digital and model predictive control.
\end{abstract}

\section{Introduction}
The method of controlled Lagrangians for stabilization of relative
equilibria (steady state motions) originated in Bloch, Leonard, and
Marsden \cite{BLM1} and was then developed in Auckly~\cite{A}, Bloch,
Leonard, and Marsden \cite{BLM2, BLM4, BLM3}, Bloch, Chang, Leonard,
and Marsden \cite{BCLM}, and Hamberg~\cite{H1, H2}. A similar approach
for Hamiltonian controlled systems was introduced and further studied
in the work of Blankenstein, Ortega, van der Schaft, Maschke, Spong,
and their collaborators (see, e.g., \cite{MaOrSc2000} and related
references). The two methods were shown to be equivalent in
\cite{ChBlLeMa2002} and a nonholonomic version was developed in
\cite{ZBM3,ZBM2002}, and \cite{Bl2003}.

According to the method of controlled Lagrangians, the original
controlled system is represented as a new, uncontrolled Lagrangian
system for a \emph{controlled Lagrangian}, a modification of the original Lagrangian. The controlled Lagrangian is designed so that its associated energy has a maximum or minimum at the (relative) equilibrium to be stabilized. The time-invariant feedback control law is obtained by requiring that  the new and old systems of equations of motion are equivalent. To obtain asymptotic stabilization, dissipation-emulating terms are added to the control input.

The method of controlled Lagrangians for discrete mechanical systems was introduced in Bloch, Leok, Marsden, and Zenkov \cite{BlLeMaZe2005}. In the present paper this formalism is further developed to include \emph{potential shaping} which is used for complete state-space stabilization of equilibria. This study is motivated by the recent development of structure-preserving algorithms for numerical simulation of controlled systems. In particular, as the closed loop dynamics of a controlled Lagrangian system is itself Lagrangian, it is natural to adopt a variational discretization that exhibits good long-time numerical stability.
  
We carry out the matching procedure explicitly for discrete systems with two degrees of freedom and prove that we can asymptotically stabilize the equilibria  of interest. The theoretical analysis is validated by simulating the discrete cart-pendulum system on an incline. When dissipation is added, the inverted pendulum configuration is asymptotically stabilized, as predicted.
We then use the discrete controlled dynamics to construct a real-time
model predictive controller with piecewise constant control
inputs. This serves to illustrate how discrete mechanics can be
naturally applied to yield digital controllers for mechanical systems.

The paper is organized as follows: In Sections \ref{discrete_mech.sec} and \ref{matching.sec} we review discrete mechanics and the method of controlled Lagrangians for stabilization of equilibria of mechanical systems. The discrete version of the potential shaping procedure and related stability analysis are discussed in Sections~\ref{discrete_matching.sec} and \ref{stabilization.sec}. The theory is illustrated with the discrete cart-pendulum system. Simulations and the construction of the digital controller are presented in Sections \ref{simulations.sec} and \ref{digital.sec}.

In a future publication we intend to treat discrete systems with
nonabelian symmetries as well as systems with nonholonomic
constraints.

\section{An Overview of Discrete Mechanics}
\label{discrete_mech.sec}
A discrete analogue of Lagrangian mechanics can be obtained by
considering a discretization of the Hamilton principle; this approach
underlies the construction of variational integrators. See Marsden and
West~\cite{MaWe2001}, and references therein, for a more detailed
discussion of discrete mechanics.

A key notion is that of the {\em discrete Lagrangian}, which is a map
$L^d: Q\times Q \rightarrow \mathbb{R}$ that approximates
the action integral along an exact solution of the Euler--Lagrange
equations joining the configurations $q_k, q_{k+1} \in Q$,
\begin{equation}
\label{exact_ld}
L^d(q_k,q_{k+1})\approx \ext_{q\in\mathcal{C}([0,h],Q)} \int_0^h
L(q,\dot q)\,dt,
\end{equation}
where $\mathcal{C}([0,h],Q)$ is the space of curves
$q:[0,h]\rightarrow Q$ with $q(0)=q_k$, $q(h)=q_{k+1}$, and $\ext$
denotes extremum. 

In the discrete setting, the action integral of
Lagrangian mechanics is replaced by an action sum
\begin{equation*}
       S ^d ( q _0, q _1, \dots, q _N) = \sum_{k=0}^{N-1}L^d (q_k, q_{k+1}),
\end{equation*}
where $q_k\in Q$, $ k = 0, 1, \dots, N $, 
is a finite sequence in the configuration space. 
The equations are obtained by the discrete Hamilton 
principle, which extremizes the discrete action given
fixed endpoints $q_0$ and $q_N$. Taking the extremum over $q_1,\dots
,q_{N-1}$ gives the {\em discrete Euler--Lagrange equations}
\begin{equation*}
         \label{DEL}
         D_1 L^d (q_k, q_{k+1}) + D_2 L^d (q_{k-1}, q_k) = 0,
\end{equation*}
for $ k = 1,\dots ,N-1$. This implicitly defines the update map
$\Phi:Q\times Q\rightarrow Q\times Q$, where
$\Phi(q_{k-1},q_k)=(q_k,q_{k+1})$ and $ Q \times Q $ replaces the phase space $ TQ $ of Lagrangian mechanics. 

Since we are concerned with control, we need to consider the effect of
external forces on Lagrangian systems. In the context of discrete
mechanics, this is addressed by introducing the {\em discrete
  Lagrange--d'Alembert principle} (see Kane, Marsden, Ortiz, and 
West~\cite{KaMaOrWe2000}), which states that
\[
\delta\sum_{k=0}^{n-1}L^{d}\left(  q_{k},q_{k+1}\right)  +\sum_{k=0}%
^{n-1}F^{d}\left( q_{k},q_{k+1}\right) \cdot\left( \delta q_{k},\delta
     q_{k+1}\right) =0
\]
for all variations $\bm{\delta q}$ of $\bm{q}$ that vanish at
the endpoints. Here, $\bm{q}$ denotes the vector of positions
$(q_0,q_1,\ldots,q_N)$, and $\bm{\delta q}
=(\delta q_0, \delta
q_1,\ldots, \delta q_N)$, where 
$\delta q_k\in T_{q_k} Q$.
The discrete one-form $F^{d}$ on $Q\times Q$ approximates the impulse
integral between the points $q_k$ and $q_{k+1}$, just as the discrete
Lagrangian $L^d$ approximates the action integral. We define 
the maps
$F^{d}_{1},F^{d}_{2}:Q\times Q\rightarrow T^{\ast}Q$ by the
relations%
\begin{align*}
F^{d}_{2}\left(  q_{0},q_{1}\right) 
\delta q_{1}
& := F^{d}\left(  q_{0},q_{1}\right) \cdot \left(  0,\delta
q_{1}\right),
\\
F^{d}_{1}\left(  q_{0},q_{1}\right) 
\delta q_{0}
& := F^{d}\left(  q_{0},q_{1}\right) \cdot \left(  \delta q_{0},0\right)  .
\end{align*}
The discrete Lagrange--d'Alembert principle may then be rewritten as%
\begin{multline*}
\delta\sum_{k=0}^{n-1}L^{d}\left(  q_{k},q_{k+1}\right) 
\\
+\sum_{k=0}%
^{n-1}\left[  F^{d}_{1}\left(  q_{k},q_{k+1}\right)  
\delta q_{k}%
+ F^{d}_{2}\left(  q_{k},q_{k+1}\right)  
\delta q_{k+1}\right]
= 0
\end{multline*}
for all variations $\bm{\delta q}$ of $\bm{q}$ that vanish at the
endpoints. This is equivalent to the {\em forced discrete
Euler--Lagrange equations}%
\begin{align*}
D_{1}L^{d}\left( q_{k}%
,q_{k+1}\right) &+ D_{2}L^{d}\left(  q_{k-1},q_{k}\right) 
\\
&+ F^{d}_{1}\left(  q_{k},q_{k+1}\right)  + F^{d}_{2}\left(
q_{k-1},q_{k}\right)  =0.
\end{align*}

\section{Matching and Controlled Lagrangians}
\label{matching.sec}
In the controlled Lagrangian approach, one considers a mechanical
system with an uncontrolled (free) Lagrangian equal to kinetic
energy minus potential energy. In the simplest
setting, we modify the kinetic energy to
produce a new controlled Lagrangian which describes the
dynamics of the controlled closed-loop system. The method
is extended by the incorporation of potential shaping in \cite{BCLM}. 

Suppose our system has configuration space $Q$ and a Lie group $G$
acts freely and properly on $Q$. It is useful to keep in mind the case
in which $Q = S \times G$ with $G$ acting only on the second factor by
the left group multiplication.
For example, for the inverted planar pendulum on a cart,
$Q =  S^1 \times \mathbb{R}$ with $G = \mathbb{R}$, the group
of reals under addition (corresponding to translations of the
cart).

Our goal is to control the variables lying in the {\em shape space} $
Q/G $ using controls that act directly on the variables lying in $G$.\footnote{The shape space is $S$ in the case $Q = S \times G$.}
For kinetic shaping, the controlled Lagrangian is constructed to be $G$-invariant, thus providing modified or {\em controlled} conservation laws.
In this paper, we assume that $G$ is an abelian group.

The key modification of the Lagrangian involves
changing the kinetic energy metric $g(\cdot,\cdot)$.
The tangent space to $Q$ can be split into a sum of horizontal
and vertical parts defined as follows: For each tangent vector
$v_q$ to $Q$ at a point $q \in Q$, we can write a unique
decomposition
$ 
   v _q  = {\rm Hor}\, v _q + {\rm Ver}\,v_q,
$ 
such that the vertical part is tangent to the orbits of the $G$-action
and the horizontal part is metric-orthogonal to the vertical
space, {\em i.e.}, it is uniquely defined by the identity
\begin{equation}
g(v_q, w_q) = g({\rm Hor}\, v_q, {\rm Hor}\, w_q) +
g({\rm Ver} \,v_q , {\rm Ver} \,w_q)
\label{mech.connection}
\end{equation}
with $v_q$ and $w_q$ arbitrary tangent vectors to
$Q$ at the point $ q \in Q $. This choice of horizontal space
coincides with that given by the \emph{mechanical connection};
see, for example, Marsden~\cite{Marsden}.

For the kinetic energy of our controlled Lagrangian, we use a
modified version of the right-hand side of equation
(\ref{mech.connection}).  The potential energy remains
unchanged. The modification consists of three ingredients:
\begin{enumerate}
\item a new choice of horizontal space, denoted
${\rm Hor}_{\tau}$,
\item  a change $g \rightarrow g _\sigma$ of the metric
on horizontal vectors,
\item a change $g \rightarrow g _\rho$
of the metric on vertical vectors.
\end{enumerate}

Let $\xi_Q$ denote the infinitesimal generator corresponding
to $\xi \in \mathfrak{g}$, where
$\mathfrak{g}$ is the Lie algebra of $G$ (see Marsden~\cite{Marsden}
or Marsden and Ratiu~\cite{MR}). Thus, for each $ \xi \in
\mathfrak{g} $, $ \xi _Q $ is a vector field on  the
configuration manifold $Q$ and its value at a point
$q \in Q$ is denoted $\xi _Q (q)$.

\begin{dfn} {\em Let $\tau$ be a Lie-algebra-valued
horizontal one-form on $Q$; that is, a one-form 
that annihilates vertical vectors. The {\bfi
$\bm \tau$-horizontal space} at $q \in Q$ consists of
tangent  vectors to $Q$ at $ q $ of the form ${\rm
Hor}_{\tau} v_q := {\rm Hor}\,v_q - [\tau (v)]_Q (q)$,
which also defines
$ v_q \mapsto {\rm Hor}_\tau v_q $, the {\bfi
$\bm \tau$-horizontal projection}. The {\bfi $\bm \tau$-vertical
projection operator} is defined by
${\rm Ver}_\tau  v_q
       := {\rm Ver} \, v_q + [\tau (v)]_Q (q)$.
}
\end{dfn}

\begin{dfn} {\em Given $g_\sigma, g_\rho$, and $\tau$,
the associated {\bfi controlled Lagrangian} $L_{\tau,\sigma,\rho}$ is given by a modified
kinetic minus the given potential energy, namely
\begin{multline} \label{controlledLag}
L_{\tau,\sigma,\rho}(v_q) =
{\textstyle \frac{1}{2} }
[ g_{\sigma}({\rm Hor}_{\tau} v_q, {\rm
Hor}_{\tau} v_q) \nonumber \\
+ g_{\rho}( {\rm Ver}_{\tau} v_q,{\rm
Ver}_{\tau} v_q) ]  - V(q).
\end{multline}
}
\end{dfn}

\medskip
The equations corresponding to this Lagrangian will be
our closed-loop equations. The new terms appearing in
those equations corresponding to the directly controlled
variables are interpreted as control inputs. The
modifications to the Lagrangian are chosen so that no
new terms appear in the equations corresponding to the
variables that are not directly controlled. We refer to
this process as \emph{matching}.

Once the form of the control law is derived using the controlled
Lagrangian, the closed-loop stability of an equilibrium can be
determined by energy methods, using any available freedom in
the choice of $\tau$, $g_{\sigma}$ and $g_{\rho}$.

We can extend the method of controlled
Lagrangians to the class of Lagrangian
mechanical systems with potential energy that may
break symmetry, \emph{i.e.}, we still have a symmetry group $G$ for the kinetic
energy of the system but we now have a potential energy
$V(x^\alpha,\theta^a)$ that need not be $G$-invariant
(see \cite{BCLM}). Further, we
consider a modification to the potential energy that also breaks symmetry
in the group variables. Let the potential energy 
for the controlled Lagrangian  be defined as
\[
V(x^\alpha,\theta^a) +
V_\epsilon(x^\alpha,\theta^a),
\]
where $V_\epsilon$ is the modification---to be determined---that
depends on a new real parameter $\epsilon$.

For many systems it is sufficient to use
the so-called \emph{simplified matching conditions} \cite{BCLM}. 
For potential shaping in the setting where the simplified matching 
conditions hold we take
$g_\rho = \rho g_{ab}$ where $\rho$ is a scalar
constant.
The controlled Lagrangian takes the form
\begin{align*}
& L _{\tau, \sigma, \rho, \epsilon} (v)  = L_{\tau, \sigma}(v) - V_\epsilon(x^\alpha,\theta^a)
+ \tfrac{1}{2} (\rho -1)  \times \nonumber 
\\
&\quad  g_{ab} \left( \dot{\theta}^a +  (g^{ac} g_{\alpha c}  + \tau_{\alpha}^a ) \dot{x} ^\alpha \right)
\left( \dot{\theta}^b +  (g^{bd} g_{\beta d}  
+ \tau_{\beta}^b ) \dot{x} ^\beta \right) , 
\end{align*}
where
\[
L_{\tau, \sigma}(v) = L(x^\alpha, \dot{x}^\beta, \theta ^a, \dot{\theta}^a
+
          \tau_{\alpha}^a \dot{x}^\alpha ) + \tfrac{1}{2}
           \sigma g_{ab} \tau_{\alpha}^a
            \tau_{\beta}^b \dot{x}^\alpha
             \dot{x}^\beta.
\]
This has sufficient generality to handle many examples of interest. 

\begin{figure}[h,t]
\begin{center}
\begin{overpic}[width=0.32\textwidth]
{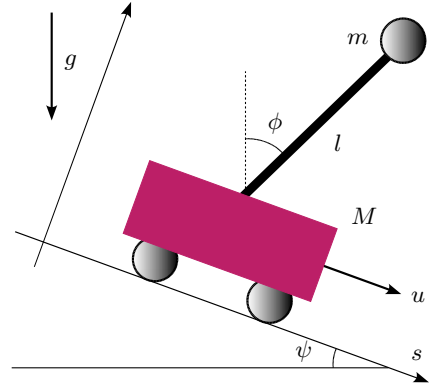}
\put(94,6){\small $s$}
\put(94,19,8){\small $u$}
\put(79,80){\small $m$}
\put(76,55){\small $l$}
\put(80,38){\small $M$}
\put(13,75){\small $g$}
\put(60.5,60){\small $\phi$}
\put(67,6.76){\small $\psi$}
\end{overpic}
\caption{\footnotesize
The pendulum on a cart going down an inclined plane under gravity. The control force is in the direction $s$, the overall motion of the cart.}
\label{cart.figure}
\end{center}
\end{figure}

A basic example treated in earlier papers in the smooth setting is
the {\em pendulum on a cart}.  Let $s$ denote the position of the cart
on the $s$-axis, $\phi$ denote the angle of the pendulum with
the upright vertical, and $\psi$ denote the elevation angle of the incline, as in Figure \ref{cart.figure}.
The configuration space for this system is $Q = S \times G = S^1
\times \mathbb{R}$, with the first factor being the pendulum angle
$\phi$ and the second factor being the cart position $s$. The
velocity phase space, $TQ$, has coordinates $(\phi, s, \dot{\phi},
\dot{s})$.  The length of the pendulum is $l$, the mass of the
pendulum is $m$ and that of the cart is $M$.

The symmetry group $G$ of the kinetic energy of the pendulum-cart system is that of translation in the $s$ variable, so $G = \mathbb{R}$. 

\section{Discrete Potential Shaping}
\label{discrete_matching.sec}
For simplicity, we consider systems with one shape and one group degree of freedom. We further assume that the configuration space $Q$ is the direct product of a one-dimensional shape space $S$ and a one-dimensional Lie group $G$. The continuous-time Lagrangian $ L : TQ \to \mathbb{R} $ and the form $\tau$ are 
\begin{equation}\label{c_lagr.eqn}
L( \phi, s, \dot \phi, \dot{s} ) 
=
{\textstyle \frac12} (
\alpha \dot \phi ^2 + 2 \beta (\phi) \dot \phi \dot{s} 
+ \gamma \dot{s} ^2 ) - V _1 (\phi) - V _2 (s)
\end{equation} 
and 
\[
\tau = k (\phi) \dot \phi \quad \text{with} \quad k (\phi) = - \frac{\beta (\phi)}{\sigma \gamma}.
\]
This Lagrangian \eqref{c_lagr.eqn} satisfies the simplified matching conditions of \cite{BCLM}. 

The continuous-time controlled Lagrangian $ L _{\tau,\sigma,\rho,\epsilon }:TQ\to \mathbb{R} $ becomes 
\begin{multline}\label{c_controlled_lagr.eqn}
\!\!\!\!
L _{\tau,\sigma,\rho,\epsilon } 
( \phi, s, \dot \phi, \dot{s} ) = 
L( \phi, s, \dot \phi, \dot{s} + k (\phi) \dot \phi) 
+
{\textstyle \frac12} \sigma \gamma 
(k (\phi) \dot \phi ) ^2 
\\
+ {\textstyle \frac12}(\rho - 1)\gamma (\dot{s} + (\sigma - 1) k (\phi) \dot \phi ) ^2 
+ V _2 (s) - V _\epsilon (y) ,\!\!
\end{multline}
where 
\[
y = s - \int _{\phi _e } ^\phi \frac{1}{\gamma}\left(
\frac{1}{\sigma} - \frac{\rho -1}{\rho} 
\right) \beta (z)\, d z,
\]
the function $ V _\epsilon (y) $ is arbitrary, and $(\phi _e, s _e) $ is the equilibrium of interest.
As in Bloch, Chang, Leonard, and Marsden \cite{BCLM},   the kinetic energies in \eqref{c_lagr.eqn} and \eqref{c_controlled_lagr.eqn} are $G$-invariant. 

For the cart-pendulum system, $\alpha$, $\beta (\phi) $, $\gamma$, $ V _1 (\phi) $, and $ V _2 (s) $ are $\alpha = m l^2$, $\beta (\phi) = ml \cos (\phi - \psi)$, $\gamma = M + m$,
$V _1 (\phi) = - mgl \cos \phi $, and $ V _2 (s) = - \gamma g s \sin \psi $. Note
that $\alpha\gamma-\beta^2 (\phi) > 0$.

In discretizing the method of controlled Lagrangians, we 
combine formulae 
\eqref{exact_ld}, \eqref{c_lagr.eqn}, and \eqref{c_controlled_lagr.eqn}.
In the rest of this paper, we will adopt the notations
\begin{equation*}\label{notations.eqn}
q _{k+1/2} = \frac {q _k + q _{k+1}}{2}, \quad 
\Delta q _k = q _{k+1} - q _k.
\end{equation*} 
This allows us to construct a \emph{second-order accurate} 
discrete Lagrangian and discrete
controlled Lagrangian as
\begin{align}
\nonumber
L ^d (q _k,q  _{k+1}) &= 
h L( q _{k+1/2} , \Delta q _k /h ),
\\
\label{discrete_second_order_controlled_lagrangian}
L^d_{\tau,\sigma,\rho,\epsilon}(q_k,q_{k+1}) &= 
h L_{\tau,\sigma,\rho,\epsilon}( 
q _{k+1/2} , \Delta q _k /h 
) ,
\end{align}
with $ q _k = ( \phi _k, s _k) $. 

The discrete dynamics is governed by the equations
\begin{align} 
\label{dcp1.eqn}
\pder{L ^d (q _k,q  _{k+1})}{\phi _k} + 
\pder{L ^d (q_{k-1},q  _k)}{\phi _k} &= 0,
\\
\label{dcp2.eqn}
\pder{L ^d (q _k,q  _{k+1})}{s _k} + 
\pder{L ^d (q _{k-1},q  _k)}{s _k} &= u _k ,
\end{align}
where $ u _k $ is the control input. 

The dynamics associated with \eqref{discrete_second_order_controlled_lagrangian} is amended by the term $ w _k $ in the discrete shape equation:
\begin{align} 
\label{controlled_lagrangian_dynamics_1.eqn}
\pder{L ^d  _{ \tau,\sigma,\rho,\epsilon} (q _k,q  _{k+1})}{\phi _k} + 
\pder{L ^d  _{ \tau , \sigma,\rho,\epsilon } (q_{k-1},q  _k)}{\phi _k} &= w _k ,
\\
\label{controlled_lagrangian_dynamics_2.eqn}
\pder{L ^d  _{ \tau , \sigma,\rho,\epsilon } (q _k,q  _{k+1})}{s _k} + 
\pder{L ^d  _{ \tau , \sigma,\rho,\epsilon } (q _{k-1},q  _k)}{s _k} &= 0. 
\end{align}
This term $ w _k $ is important for matching systems \eqref{dcp1.eqn}, \eqref{dcp2.eqn} and \eqref{controlled_lagrangian_dynamics_1.eqn}, \eqref{controlled_lagrangian_dynamics_2.eqn}. \emph{The presence of the terms $w _k$ represents an interesting (but manageable) departure from the continuous theory.}
Let 
\[
J _k = \rho \gamma 
\big(
{\Delta s _k}/{h} 
- (\sigma - 1) 
k(\phi _{k+\frac12}) 
{\Delta \phi _k}/{h}
\big).
\]
The following statement is proved by a straightforward calculation:
\begin{theorem}\label{discrete_matching.thm}
{\em
The dynamics \eqref{dcp1.eqn}, \eqref{dcp2.eqn}  is equivalent to the dynamics
\eqref{controlled_lagrangian_dynamics_1.eqn},
\eqref{controlled_lagrangian_dynamics_2.eqn} if and only if $u _k$ and $w _k$ are given by
}
\begin{align}\label{discrete_u.eqn}
\nonumber
u _k &= 
- \frac{h}{2} \left[
V _2 ' (s _{k+\frac12}) + V _2 ' (s _{k-\frac12})
\right]
\\
\nonumber
& \quad \,
+ \frac{h}{2 \rho} \left[
V _\epsilon  ' (s _{k+\frac12}) + V _\epsilon ' (s _{k-\frac12})
\right]
\\
& \quad \, +
\frac{\gamma \Delta \phi _k k ( \phi _{k+1/2} ) 
- \gamma \Delta \phi _{k-1} k ( \phi _{k-1/2})}{h },
\end{align}
and
\begin{align*} 
w _k &= -
\Big(
1 - \sigma  + \frac{\sigma}{\rho}
\Big)
\Big(
k( \phi _{k+\frac12})
\Big[ - \gamma \rho J _k + 
\frac{h}{2} V ' _\epsilon (y _{k+\frac12})
\Big]
\\
& \quad \, 
+ k( \phi _{k-\frac12})
\Big[ \gamma \rho J _{k-1} + 
\frac{h}{2} V ' _\epsilon (y _{k-\frac12})
\Big]
\\
& \quad \, 
- k ' ( \phi _{k+\frac12}) J _k \Delta \phi _k 
- k ' ( \phi _{k-\frac12}) J _{k - 1} 
\Delta \phi _{k - 1}
\Big) . 
\end{align*} 
\end{theorem}
\smallskip

\emph{Remark.\ }
The terms $ w _k $ vanish when $ \beta (\phi) = \text{const} $ as they become proportional to the left-hand side of equation \eqref{controlled_lagrangian_dynamics_2.eqn}. 

\section{Stabilization of the Discrete Controlled System}
\label{stabilization.sec}
The stability analysis in this paper is done by means of an analysis of the spectrum of the linearized discrete equations. We assume that the equilibrium to be stabilized is $ (\phi _k, s _k) = (0,0) $.  

\begin{theorem}
\emph{
The equilibrium $ (\phi _k , s _k) = (0, 0) $ of equations \eqref{controlled_lagrangian_dynamics_1.eqn} and \eqref{controlled_lagrangian_dynamics_2.eqn} is spectrally stable if 
\begin{equation}\label{discrete_spectrum.eqn}
- \frac{\beta ^2 (0)}{\alpha \gamma - \beta ^2 (0)} < \sigma < 0, \quad \rho < 0, \quad \text{and} \quad V _\epsilon '' (0) < 0.\!
\end{equation}
}
\end{theorem}
\medskip
\begin{proof}
The linearized discrete equations are 
\begin{align} 
\label{linearized_1.eqn}
\pder{\mathcal L ^d  _{ \tau,\sigma,\rho,\epsilon} (q _k,q  _{k+1})}{\phi _k} + 
\pder{L ^d  _{ \tau , \sigma,\rho,\epsilon } (q_{k-1},q  _k)}{\phi _k} &= 0 ,
\\
\label{linearized_2.eqn}
\pder{\mathcal L ^d _{ \tau , \sigma,\rho,\epsilon } (q _k,q  _{k+1})}{s _k} + 
\pder{\mathcal L ^d _{ \tau , \sigma,\rho,\epsilon } (q _{k-1},q  _k)}{s _k} &= 0, 
\end{align}
where $ \mathcal{L} ^d _{ \tau , \sigma,\rho,\epsilon } (q _k,q  _{k+1}) $ is the quadratic approximation of $ L ^d _{ \tau , \sigma,\rho,\epsilon } $ at the equilibrium 
(\emph{i.e.}, $ \beta (\phi) $,  $ V _1 (\phi) $, and $ V _\epsilon (y) $ in $ L ^d _{ \tau , \sigma,\rho,\epsilon } $ are replaced by $ \beta (0) $, $ \frac12 V _1 '' (0) \phi ^2 $, and  $ \frac 12 V _\epsilon '' (0) y ^2 $, respectively). \emph{Note the absence of the term $ w _k $ in equation~\eqref{linearized_1.eqn}. }

The linearized dynamics preserves the quadratic approximation of the
discrete energy
\begin{multline}\label{quadratic_energy.eqn} 
\frac{\alpha \gamma \sigma ^2 - 
\beta (0) ^2 (\sigma-1)(\rho (\sigma-1) - \sigma )}
{2 \gamma \sigma ^2 h} 
\Delta \phi _k^2 
\\
+ \frac{\beta (0) \rho (\sigma - 1)}{\sigma h}
\Delta \phi _k \Delta s _k 
+ \frac{\gamma \rho}{2 h} \Delta s _k ^2 
\\
+ \frac{h}{2} V _1 '' (0) \phi _{k+\frac12} ^2 
+ \frac{h}{2} V _\epsilon '' (0) x _{k+\frac12} ^2 
,
\end{multline} 
where
\[ 
x  = s  + 
\left(
\frac{\rho - 1}{\rho} - \frac{1}{\sigma}
\right)
\frac{\beta (0)}{\gamma}\, \phi .
\]
Since $ V _1 '' (0) $ is negative, 
the equilibrium $ (\phi _k, s _k) = (0,0) $ of equations 
\eqref{linearized_1.eqn} and \eqref{linearized_2.eqn} is stable if 
the function \eqref{quadratic_energy.eqn} is negative-definite. The latter requirement is equivalent to
conditions \eqref{discrete_spectrum.eqn}. The spectrum of the linearized discrete dynamics in this case belongs to the unit circle. \end{proof}
\smallskip

\emph{Remarks.\ }
Spectral stability in this situation is not sufficient to conclude
nonlinear stability.
The stability conditions \eqref{discrete_spectrum.eqn} are identical to the
stability conditions of the corresponding continuous-time system.

Following \cite{BCLM}, we now modify the control input \eqref{discrete_u.eqn} by adding the
{\em discrete dissipation-emulating term} 
\begin{equation}\label{friction.eqn}
- \frac{D (\Delta y _{k-1} + \Delta y _k)}{h}
\end{equation} 
in order to achieve the asymptotic stabilization of the equilibrium $ (\phi _k, s _k) = (0,0) $. In the above, $D$ is a constant. 
The linearized discrete dynamics becomes
\begin{align} 
\nonumber 
\label{linearized_damped_1.eqn}
\pder{\mathcal L ^d  _{ \tau,\sigma,\rho,\epsilon} (q _k,q  _{k+1})}{\phi _k} + 
\pder{L ^d  _{ \tau , \sigma,\rho,\epsilon } (q_{k-1},q  _k)}{\phi _k} \qquad \quad &
\\
= - \left( 
\frac{\rho - 1}{\rho} - \frac{1}{\sigma}
\right)
\frac{\beta (0)}{\gamma} 
\frac{D (\Delta x _{k-1} + \Delta x _k)}{h} &,
\\
\nonumber
\label{linearized_damped_2.eqn}
\pder{\mathcal L ^d _{ \tau , \sigma,\rho,\epsilon } (q _k,q  _{k+1})}{s _k} + 
\pder{\mathcal L ^d _{ \tau , \sigma,\rho,\epsilon } (q _{k-1},q  _k)}{s _k} \qquad \quad &
\\
= - \frac{D (\Delta x _{k-1} + \Delta x _k)}{h} &. 
\end{align}
\begin{theorem} \label{asymtotic_stability.thm}
\emph{
The equilibrium $ (\phi _k, s _k) = (0, 0) $ of equations \eqref{linearized_damped_1.eqn} and \eqref{linearized_damped_2.eqn} is asymptotically stable if conditions \eqref{discrete_spectrum.eqn} are satisfied and $ D $ is positive.
}
\end{theorem} 
\begin{proof}
Multiplying equations \eqref{linearized_damped_1.eqn} and \eqref{linearized_damped_2.eqn} by $ (\Delta \phi _{k-1} + \Delta \phi _k )/2 $ and $ (\Delta s _{k-1} + \Delta s _k )/2 $, respectively, we obtain 
\[
E _{k, k+1 } = E _{k-1, k } + \frac{D h}{4} 
\left(
\frac{\Delta x _{k-1}}{h} + \frac{\Delta x _k}{h}
\right) ^2,
\] 
where $ E _{k, k+1} $ is the quadratic approximation of the discrete energy \eqref{quadratic_energy.eqn}. Recall that $ E _{k, k+1} $ is negative-definite. It is possible to show that, in some neighborhood of $ (\phi _k, s _k) = (0,0) $,  the quantity $ \Delta x _{k-1} + \Delta x _k \not\equiv 0 $ along a solution of equations \eqref{linearized_damped_1.eqn} and \eqref{linearized_damped_2.eqn} unless this solution is the equilibrium $ (\phi _k, s _k) = (0,0) $. Therefore, $  E _{k, k+1} $ increases along non-equilibrium solutions of \eqref{linearized_damped_1.eqn} and \eqref{linearized_damped_2.eqn}. Since equations \eqref{linearized_damped_1.eqn} and \eqref{linearized_damped_2.eqn} are linear, this is only possible if the spectrum of \eqref{linearized_damped_1.eqn} and \eqref{linearized_damped_2.eqn} is inside the open unit disk, which implies asymptotic stability of the equilibrium of both linear system \eqref{linearized_damped_1.eqn} and \eqref{linearized_damped_2.eqn} and nonlinear system \eqref{dcp1.eqn} and \eqref{dcp2.eqn} with discrete dissipation-emulating term \eqref{friction.eqn} added to $ u _k $.
\end{proof}

\section{Simulations}\label{simulations.sec}

Simulating the discrete behavior of the controlled Lagrangian system
involves viewing equations \eqref{dcp1.eqn} and \eqref{controlled_lagrangian_dynamics_2.eqn} as an implict update map $\Phi:(q_{k-2},
q_{k-1})\mapsto(q_{k-1},q_k)$. This presupposes that the initial
conditions are given in the form $(q_0,q_1)$; however it is generally
preferable to specify the initial conditions as $(q_0,\dot q_0)$. 
This is achieved by solving the 
boundary condition
\[
\frac{\partial L}{\partial \dot q}(q_0,\dot q_0) + D_1 L^d(q_0,q_1) +
F^d_1(q_0,q_1)
= 0
\]
for $q_1$. Once the initial conditions are expressed in the form $(q_0,q_1)$, the discrete evolution can be obtained using the implicit update
map $\Phi$.

\begin{figure}[h]
\hspace{-1.15em}
\begin{overpic}
[scale=.123]
{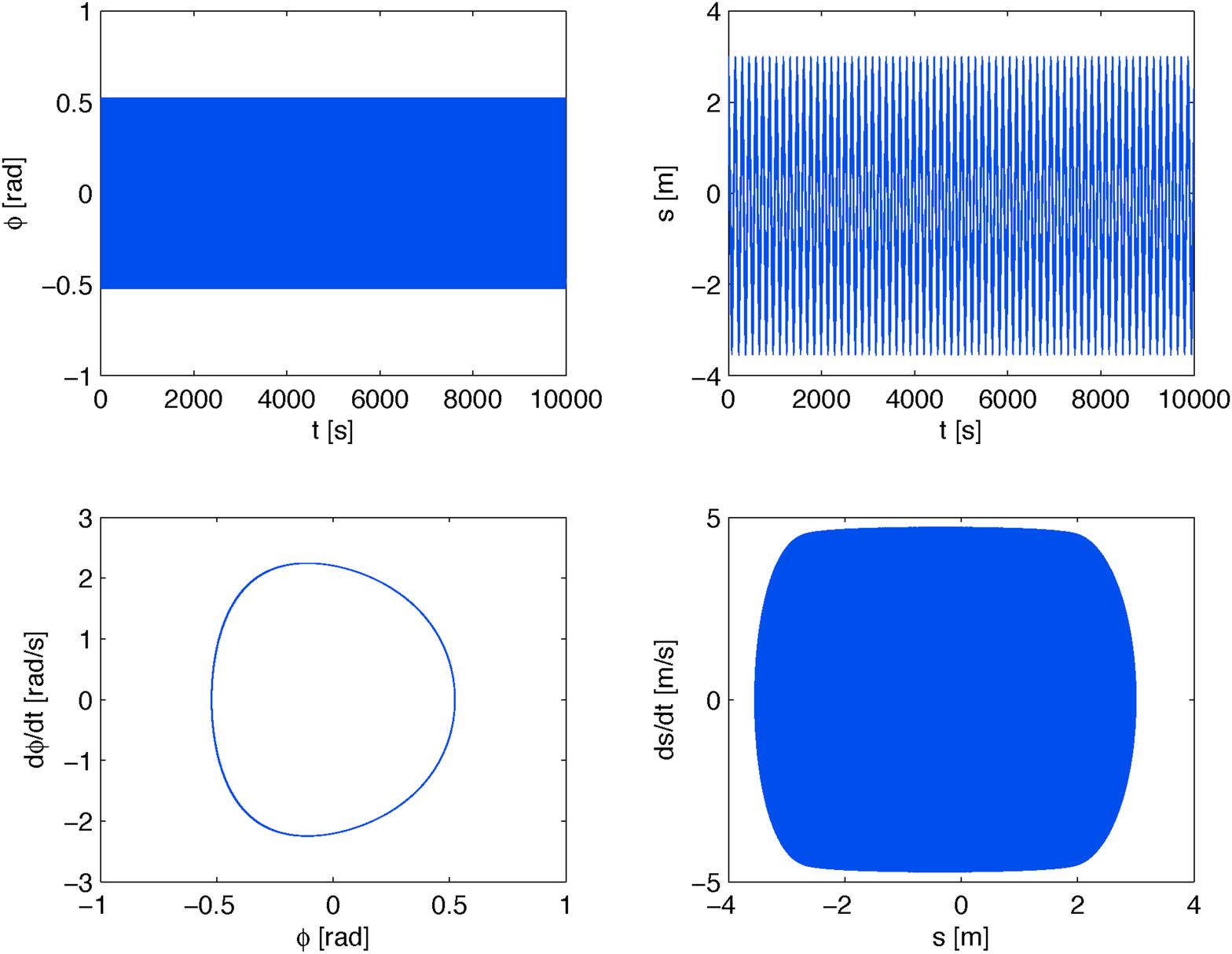}
\put(25.8,10.55)
	{\begin{rotate}{0}
	{\includegraphics[width=.015\textwidth]
	{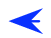}}
	\end{rotate}
	}
\put(29,31.67)
	{\begin{rotate}{180}
	{\includegraphics[width=.015\textwidth]
	{Figures/arrow}}
	\end{rotate}
	}
\end{overpic} 
\caption{Discrete controlled dynamics without dissipation. The discrete controlled system stabilizes the motion about the equilibrium; since there is no dissipation, the oscillations are sustained. }
\label{nodiss}
\end{figure}

In Figure~\ref{nodiss}, we present a MATLAB simulation of discrete controlled dynamics of the cart-pendulum system in the absence of dissipation. 

Here, $h=0.05\,\textrm{sec}$, $m=0.14\,\textrm{kg}$, $M=0.44\,\textrm{kg}$, $l=0.215\,\textrm{m}$, and $\psi=\frac{\pi}{9}\,\textrm{radians}$. Our goal is to regulate the cart at $s=0$ and the pendulum at $\phi=0$. The control gains are chosen to be $\kappa=20$, $\rho=-0.02$, and $\epsilon=0.00001$. It~is worth noting that the discrete dynamics remain bounded near the desired equilibrium, and this behavior persists even for significantly longer simulation runs involving $10^6$ time-steps. The exceptional stability of the discrete controlled trajectory can presumably be understood in terms of the bounded energy oscillations characteristic of symplectic and variational integrators.

\begin{figure}[htbp]
\hspace{-1.15em}
\begin{overpic}
[scale=.51]
{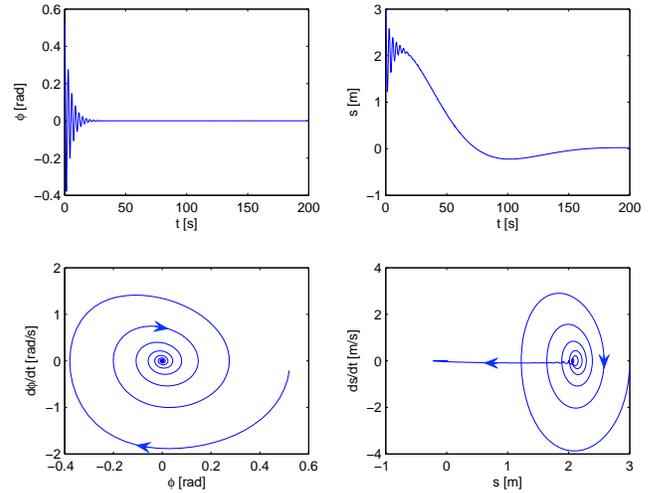}
\put(28,26.45)
	{\begin{rotate}{170}
	{\includegraphics[width=.015\textwidth]
	{Figures/arrow}}
	\end{rotate}
	}
\put(22.1,8.64)
	{\begin{rotate}{-10}
	{\includegraphics[width=.015\textwidth]
	{Figures/arrow}}
	\end{rotate}
	}
\put(70,19.85)
	{\begin{rotate}{-1}
	{\includegraphics[width=.015\textwidth]
	{Figures/arrow}}
	\end{rotate}
	}
\put(88.13,19.5)
	{\begin{rotate}{90}
	{\includegraphics[width=.015\textwidth]
	{Figures/arrow}}
	\end{rotate}
	}
\end{overpic}
\caption{Discrete controlled dynamics with dissipation. Here the oscillations die out and the cart converges to the desired point $s = 0$.}
\label{diss}
\end{figure}

When dissipation is added, we obtain an asymptotically stabilizing control law, as illustrated in Figure~\ref{diss}. This is consistent with the stability analysis of Section \ref{stabilization.sec}.

\section{Model Predictive Controller}\label{digital.sec}

We now explore the use of the forced discrete Euler--Lagrange
equations as the model in a real-time model
predictive controller, with piecewise constant control
forces. Algorithm 1 below describes the details of the procedure.
\begin{algorithm}
\caption{\textsc{Digital Controller} ( $q(\,\cdot\,), T_f, h$ )}
\begin{algorithmic}
\STATE $q_0\leftarrow$ \textbf{sense} $q(0)$
\STATE $q_1\leftarrow$ \textbf{sense} $q(h)$
\STATE $\bar q_2 \leftarrow$ \textbf{solve} $D_2 L^d(q_0, q_1) + D_1
L^d(q_1,\bar q_2) =0$ 
\STATE $\bar q_3 \leftarrow$ \textbf{solve} $D_2 L^d(q_1, \bar q_2) +
D_1 L^d(\bar q_2,\bar q_3) + F^d_1(\bar q_2, \bar q_3) =0$ 
\STATE $u_{2+1/2} \leftarrow u\left(\frac{\bar q_2 + \bar q_3}{2},\frac{\bar q_3-\bar q_2}{h}\right)$  
\STATE \textbf{actuate} $u=u_{2+1/2}$ for $t\in[2h,3h]$
\STATE $q_2 \leftarrow$ \textbf{sense} $q(2h)$
\STATE $\bar q_3 \leftarrow$ \textbf{solve} $D_2 L^d(q_1, q_2) + D_1
L^d(q_2,\bar q_3) + F_1^d(q_2, \bar q_3) =0$ 
\STATE $\bar q_4 \leftarrow$ \textbf{solve} $D_2 L^d(q_2, \bar q_3) +
D_1 L^d(\bar q_3,\bar q_4)$
\\
\hspace*{1in}$+ F_2^d(q_2, \bar q_3) + F_1^d(\bar q_3, \bar q_4) =0$
\STATE $u_{3+1/2} \leftarrow u\left(\frac{\bar q_3 + \bar q_4}{2},\frac{\bar q_4-\bar q_3}{h}\right)$
\STATE \textbf{actuate} $u=u_{3+1/2}$ for $t\in[3h,4h]$
\FOR{$k=4$ to $(T_f/h -1)$}
\STATE $q_{k-1} \leftarrow$ \textbf{sense} $q((k-1)h)$
\STATE $\bar q_k \leftarrow$ \textbf{solve} $D_2 L^d(q_{k-2}, q_{k-1})
+ D_1 L^d( q_{k-1},\bar q_k)$
\\
\hspace*{1in}$+ F_2^d(q_{k-2}, q_{k-1}) + F_1^d( q_{k-1}, \bar q_k) =0$
\STATE $\bar q_{k+1} \leftarrow$ \textbf{solve} $D_2 L^d(q_{k-1}, \bar
q_k) + D_1 L^d(\bar q_k,\bar q_{k+1})$
\\
\hspace*{1in}$+ F_2^d(q_{k-1}, \bar q_k) + F_1^d(\bar q_k, \bar q_{k+1}) =0$
\STATE $u_{k+1/2} \leftarrow u\left(\frac{\bar q_k + \bar q_{k+1}}{2},\frac{\bar q_{k+1}-\bar q_k}{h}\right)$
\STATE \textbf{actuate} $u=u_{k+1/2}$ for $t\in[kh,(k+1)h]$
\ENDFOR
\end{algorithmic}
\end{algorithm}

The digital controller uses the position information it senses for
$t=-2h, -h$ to estimate the positions at $t=0,h$ during the time
interval $t=[-h,0]$. This allows it to compute a symmetric finite
difference approximation to the continuous control force $u(\phi,s,\dot\phi,\dot s)$ at $t=h/2$ using the
approximation
\begin{align*}
u_{1/2} &= u\left(\frac{\bar\phi_0+\bar\phi_1}{2},\frac{\bar s_0+\bar s_1}{2}, \frac{\bar\phi_1-\bar\phi_0}{h},\frac{\bar s_1-\bar s_0}{h}\right),
\end{align*}
where the overbar indicates that the position variable is being
estimated by the numerical model. This control is then applied as a
constant control input for the time interval $[0,h]$. This algorithm can be implemented in real-time if the two
forward solves can be computed within the time interval $h$.

The initialization of the discrete controller is somewhat involved,
since the system is unforced during the time interval $[0,2h]$ while
the controller senses the initial states, and computes the appropriate
control forces. 

\begin{figure}[thbp]
\hspace{-1.15em}
\begin{overpic}
[scale=.51]
{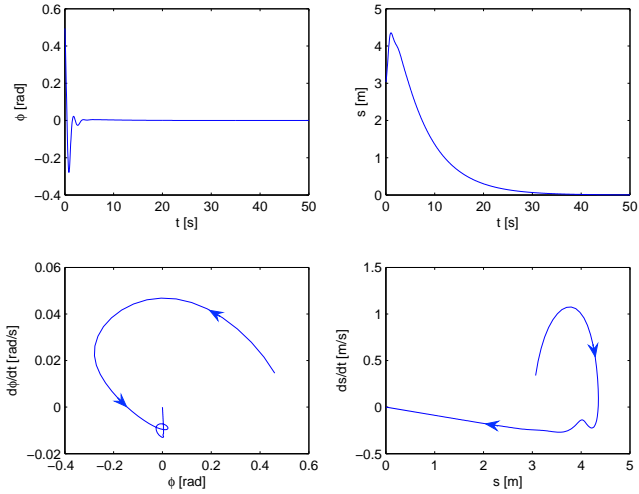} 
\put(22.7,15)
	{\begin{rotate}{135.3}
	{\includegraphics[width=.015\textwidth]
	{Figures/arrow}}
	\end{rotate}
	}
\put(31.8,27.48)
	{\begin{rotate}{-30.5}
	{\includegraphics[width=.015\textwidth]
	{Figures/arrow}}
	\end{rotate}
	}
\put(70,11.52)
	{\begin{rotate}{-9.3}
	{\includegraphics[width=.015\textwidth]
	{Figures/arrow}}
	\end{rotate}
	}
\put(87,21.15)
	{\begin{rotate}{100}
	{\includegraphics[width=.015\textwidth]
	{Figures/arrow}}
	\end{rotate}
	}
\end{overpic}

\caption{The discrete real-time piecewise constant model predictive controller stabilizes $\phi$ and $s$ to zero.}
\label{digital_control}
\end{figure}

The numerical simulation of the digital controller is shown in
Figure~\ref{digital_control}. We see that the system is \nolinebreak
asymptotically \nolinebreak stabilized \nolinebreak in both the $\phi$ and $s$ variables.

\section{Conclusions}
In this paper we have introduced potential shaping techniques 
for discrete systems and have shown that these lead to an effective
numerical implementation for stabilization in the case of the 
discrete cart-pendulum model. 
The method in this paper is related to other discrete methods in control that have a long history; recent papers that use discrete mechanics 
in the context of optimal control and celestial navigation 
are~\cite{GuBl2005}, \cite{JuMaOb2005}, and \cite{SaShMcBl2005}. The full theory of discrete controlled Lagrangians will be developed in a forthcoming paper. 

\vspace{3.9em}

\section{Acknowledgments}
The authors would like to thank the reviewers for helpful remarks.
The research of AMB was supported by NSF grants DMS-0305837, DMS-0604307, and CMS-0408542. The research of ML was partially supported by NSF grant DMS-0504747 and a University of Michigan Rackham faculty grant. The research of JEM was partially supported by AFOSR Contract FA9550-05-1-0343. The research of DVZ was partially supported by NSF grants DMS-0306017 and DMS-0604108.

\end{document}